\documentclass[12pt,a4paper,twoside]{amsart}
\usepackage{amsmath}
\usepackage{amssymb}
\usepackage{amsfonts}

\date{}

\def\SL{\rm{SL}}

\def\Ker{\text{\rm Ker}}
\def\Aut{\text{Aut}}

\def\vp{\varphi}

\def\bbr{\mathbb{F}}

\theoremstyle{plain}
\newtheorem{thm}{Theorem}
\newtheorem*{thm*}{Theorem}

\theoremstyle{definition} 

\newtheorem*{definition*}{Definition}
\newtheorem*{thm1*}{Theorem A1}
\newtheorem*{thm2*}{Theorem A2}
\newtheorem{cor}[thm]
{Corollary}

\newtheorem*{claim*}{Claim}
\newtheorem*{rem*}{Remarks}

\def\bbr{\mathbb{R}}

\def\be{\begin{equation}}
\def\ee{\end{equation}}

\def\a{\alpha}
\def\vare{\varepsilon}
\def\ol{\overline}

\theoremstyle{remark}  
\begin{document}

\title{Images of word maps \break in finite simple groups}
\author[Alexander Lubotzky]{Alexander Lubotzky$^*$}
\address{Institute of Mathematics, Hebrew University\\
Jerusalem 91904, ISRAEL\\
alex.lubotzky@mail.huji.ac.il}
\thanks{*This work was supported by ERC, NSF and ISF}
\maketitle

\begin{abstract}
In response to questions by Kassabov, Nikolov and Shalev, we show that a given subset $A$ of a finite simple group $G$ is the image of some word map $w : G\times G\to G $ if and only if (i) $A$ contains the identity and (ii) $A$ is invariant under $\Aut (G)$.
\end{abstract}

\section{Introduction}

Let $w$ be a word in the free group $F_k$ of rank $k$.  For every group $G$, $w$ defines a ``word map" $w : G^k\to G$.  Let $w(G)$ denote the image of this map.  In recent years there has been  great interest in $w(G)$, especially when $G$ is a finite simple group.  For example, the  Ore conjecture \ \ \ \ \ \  settled in [LOST] says that for $w = xy x^{-1} y^{-1}, w(G) = G$ for every non-abelian finite simple group.  Similar results have been proved for a few  other words.  Another direction of research is: what is the width of $G$ with respect to $w(G)$?  In this respect, the most remarkable result [LST] is that for every given $w$ the width is two (i.e., every element of $G$ is a product of two elements of $w(G)$) if $G$ is large enough (depending on $w$).   For more, see the above mentioned papers, [Se], [W] and the references therein.

In these results, the word is fixed and $G$ is changed. Several recent papers have been devoted to the dual question: Given $G$ what kind of subsets can appear as $w(G)$ for some $w$.  For example, Kassabov and Nikolov [KN] showed that in $A_n, n \ge 7$, the set consisting of $e$ and all the 3-cycles is $w(A_n)$ for some $w$ (depending on $n$), which shows, in particular, that the width of $G$ with respect to $w(G)$ can be arbitrary large.  They also showed  for $G = \SL_n(q)$, that there exists a word $w$ with $w(G)$ equal the identity and all transvections. (A similar result for some sequences of $n$ and $q$ has been proved earlier by Guralnick and Tiep.)  Levy [L] showed that for $q=2^{2^n}$, $n\ge 2$ and $G = \SL_2 (q)$, there exists $w \in F_2$ for which $w(G)$ consists of the identity plus the union of four conjugacy classes of elements of order 17.  See also [JLO] for related  results. Some questions have been asked about the possible subsets of $G$, a finite simple group, to be equal to $w(G)$ for some $w$.

There are two clear necessary conditions for a subset $A \subseteq G$ to be equal to $w(G)$:

\begin{enumerate}
\item $e\in A$ (since $w(e, \dots, e) = e)$.

\item For every $\a \in \Aut (G), \a (A) = A $ (since $\a (w(g_1,\dots, g_k)) =w(\a(g_1), \dots, \a (g_k))$.

\end{enumerate}

During the conference ``Words and Growth" (Jerusalem June 2012), Shalev asked whether this could also be sufficient. The goal of this note is to answer this question affirmatively.

\begin{thm}  Let $G$ be a finite simple group and $A$ a subset of $G$ such that $e \in A$ and for every $\a \in \Aut (G), \a (A) = A$.  Then there exists a word $w\in F_2$, s.t. $w(G) = A$.
\end{thm}
The proof is fairly elementary but we make use, along the way, of a result of Guralnick and Kantor [GK, Corollary, p.~745] asserting that  for every finite simple group $G$ and for every $e \neq a \in G$ there exists $b\in G$ such that $G = \langle a, b\rangle$, i.e. $G$ is generated by $a$ and $b$.  The proof of this result requires the classification of the finite simple groups.  Thus our result also depends on the classification.  While it seems impossible at this stage to prove the Guralnick-Kantor result without the classification (in fact, without it, it is not even known that every finite simple group is generated by two element) it might be that our result has a classification free proof.

The method of proof has an interesting Corollary:
\begin{cor} For every finite simple group $G$, there is $w(x,y)\in F(x,y)$, the free group on two generators with the following property:
\[ \mathrm{For\ every\ } \; (a, b) \in G\times G, \; \; w(a,b) \neq e\, \mathrm{\ if\ and\ only\ if\ } \langle a, b\rangle = G .\]
\end{cor}

So $w$ can ``test" whether two elements generate $G$.  Unfortunately our method of proof while ``effective" (see a remark at the end of \S 2) does not really give a useful description of the word $w$ of Theorem 1 or Corollary 2.  So after all, the methods of [KN] and [L] have some advantage in spite of proving only very special cases of the theorem.

\medskip

\noindent{Acknowledgment.}  The author is grateful to R. Guralnick, M. Kassabov,  N. Nikolov and A. Shalev for some interesting discussions during and after the above noted conference.

\section{Proof of the theorem}

Let $G$ be a finite simple group and $A \subset G$ with $e \in A$ and $\a (A) = A$ for every $\a \in \Aut (G)$.  If $G$ is abelian then clearly either $A = \{ e\}$ or $A=G$ and in both cases the result is trivial, so from now on we assume $G$ is non abelian.

Let $\{ (a_i, b_i) | i = 1, \dots, |G|^2\}$ be the set of all ordered pairs of elements of $G$ such that for $i= 1,\dots, \ell$, the pair $\{ a_i, b_i\}$ generates $G$ i.e., $G = \langle a_i, b_i\rangle$, while  for $i = \ell +1,\dots, |G|^2$, it does not.  Every $i$ gives rise to a unique homomorphism $\varphi_i:F_2 \to G$ defined by $\vp_i(x) = a_i, \vp_i  (y) = b_i$ when $F = F_2$ is the free group  on $x$ and $y$. Denote $N_i = \Ker \vp_i$ and let $M = \bigcap^{|G|^2}_{i = 1} N_i$ and $N = \bigcap^\ell_{i = 1} N_i$.
So for $i = 1,\dots, \ell,\; \vp_i $ is an epimorphism and $F/N_i \simeq G$, while for $i > \ell, F/N_i$ is isomorphic to a proper subgroup of $G$.  Let $\vp_M = (\vp_1,\dots, \vp_{|G|^2})$ be the diagonal map from $F $ to $G^{|G|^2}$ and $H = \vp_M (F_2)$.  We want to describe the structure of $H$ as a subgroup of $G^{|G|^2}$. Write $G^{|G|^2}$ as $E\times D$ where $E$ is the product of the first $\ell $ copies of $G$ (the ones corresponding to epimorphisms to $G$) and $D$ is the product of all the others.  Let $K$ be the kernel of the projection from $H$ to $E$ and $E'$ its image there, and let $L$ be the kernel of the projection from $H$ to $D$ and $D'$ will denote the image.  The group $K$ is a subgroup of $D$, whose projection to every single copy of $G$ in $D$, is a proper subgroup of $G$.  So $K$ has no Jordan-H\"older factor isomorphic to $G$. On the other hand $E' \cong H/K$ which is the projection of $H$ to $E$ is a subdirect product of $G^\ell$ such that its projection to every single copy of $G$ is onto.  Hence $E'$  is isomorphic (since $G$ is finite and simple) to $G^{r'}$ for some $r' \le \ell$.  We can determine precisely what is $r'$.  Let us postpone this computation for a moment, but observe first that $L$, the  kernel of the projection from $H$ to $D$, must be equal to the projection $E'$ of $H$ to $E$ since it is a subgroup of this projection and both are isomorphic to $G^{r'}$, since all the $r'$ Jordan-H\"older factors of $H$ should appear in $L$, as $H/L$ has no Jordan-H\"older factor isomorphic to $G$. Thus $H=E' \times D'$.  Note also that an element $\overline u = (u_1,\dots, u_\ell)$ is in $E'$ if and only if whenever $\a \circ \vp_i = \vp_j$ for some $1 \le i, j \le \ell$ and $\a \in \Aut (G)$, $\a (u_i) = u_j$.

Let us now calculate $r'$:

The group $\Aut (G)$ acts on the pairs $\{ (a_i, b_i)\}$ and similarly on the homomorphisms $\{ \vp_i\}$, preserving the first $\ell$ of them (the epimorphisms).  The action on these epimorphisms is free: indeed, if $\a \in \Aut (G) $ and $\a \circ \vp_i = \vp_i$ (or equivalently: $(\a (a_i), \a (b_i)) = (a_i, b_i)$) then $\a$ is the identity automorphism of $G$.  Thus, the first $\ell$ homomorphisms i.e., the $\ell$ epimorphisms, form $r = \frac{\ell}{|\Aut (G)|}$ orbits.  Now $F_2 /N \simeq H/K$ is the maximal quotient of $F_2$ which is isomorphic to a direct power $G^{r'}$ of $G$.  By a result of P.~Hall - see [KL, Corollary 7], $r' = r = \frac{\ell}{|\Aut (G)|}.$

In summary, the group $H = \vp_M (F_2)$ is a direct product $H= E' \times D'$ where $E'$ is a subgroup of $E=G^\ell$  isomorphic to $G^r$, embedded ``diagonally" in $G^\ell $ twisted by $\Aut (G)$. The other part $D'$ is a subgroup of $D = G^{|G|^2-\ell}$ whose structure is less clear, but $D'$ has no Jordan-H\"older factor isomorphic to $G$.

Let us now look at the set $A' = A \setminus \{ e \}$.  This set is a union of orbits of $\Aut (G)$ acting on $G$.  We first observe that the number of orbits is less or equal $r$.  Indeed, by the  Guralnick-Kantor result mentioned in the introduction [GK, Corollary, p.~745] every  $a \in A'$ is part of a two-element set of generators, so there exists at least one $b \in G$ such that $\langle a, b\rangle = G$ and so there exists $1\le i \le \ell$ s.t. $(a_i, b_i) = (a, b)$.  The orbit of $a$ in $G$, i.e. $\{ \a (a) | \a \in \Aut (G)\}$, gives an orbit of pairs (equivalently, of epimorphisms) $\{ ( \a (a_i), \a (b_i)) | \a \in \Aut (G)\}$.  In general, there may be more than one orbit of epimorphisms corresponding to $a$ as there may be $b$ and $b'$ such that $\langle a, b\rangle = \langle a, b'\rangle = G$ while $(a, b)$ and $(a, b')$ are not on the same $\Aut(G)$ orbit.

Let us now define the following element $\ol z = (z_i)^{|G|^2}_{i = 1} $ of $G^{|G|^2} = E \times D$:

\[ z_i = \begin{cases} a_i &\text{if} \; i \le \ell \; \text{and} \; a_i \in A'\\
e &\text{otherwise}\end{cases}.\]

We first claim that $\ol z$ is in $H=E' \times D'$.  Clearly its projection to $D$ is the identity, so we just need to check that its projection to $E$ is in $E'$.  But this is clear as by its definition, $\ol z$ is on the ``twisted diagonal" group defining $E'$ as a subgroup of $E$, and $A'$ is $\Aut (G)$-invariant.  We also observe that all the elements of $A'$ show up as some coordinates of $\ol z$.  This follows from the Guralnick-Kantor result which ensures that every $a \in A'$ has a mate $b$, with $(a,b) = (a_i, b_i)$ for some $1\le i \le \ell$.

The description of $H$ above shows that $\ol z$ is an element of $H = \vp_M (F_2) \in G^{|G|^2}$.  Spelling out the meaning of this, we see that there exists $w\in F_2$ such that for every $i = 1,\dots, |G|^2, \; \vp_i (w) = z_i$.  But $\vp_i(w) = w (a_i, b_i)$.  This means that $w$ is a word in $F_2$, with $w(G) = \{ z_i\}^{|G|^2}_{i = 1} = \{ e\} \cup A' = A$ and the theorem is proved.

The proof shows that \[w(a, b) = \begin{cases} a\in A' &\text{and} \; \langle a, b\rangle = G\\
e &\text{otherwise}\end{cases}.\]  Applying this for the set $A=G$, we deduce Corollary 2.


\bigskip

\noindent {\bf Remarks \ 1)} \ The proof actually shows that if $f : G\times G \to G$ is any function which is $\Aut (G)$-invariant (i.e., for every $\a \in \Aut \, G$, and every $a, b \in G, \; f \big(\a (a), \a(b)) = \a (f(a, b)\big)$) and $f(a, b) = e$ if $\langle a, b\rangle \neq G$, then there exists $w \in F_2$ such that $f(a, b) = w(a, b)$ for all $a, b \in G$. One can take for example $$f (a, b) = \begin{cases} a b a^{-1} b^{-1} &\text{if} \; \; \langle a, b\rangle = G \\
e &\text{otherwise}\end{cases}.$$

\medskip

\noindent {\bf 2)} \  Recall the well known result that almost every pair $(a, b) \in G\times G$ generates $G$.  From this and Remark 1, one can deduce that for every $\vare > 0$, if $G$ is large enough and $p$ a probability function on $G$ which is $\Aut (G)$-invariant (i.e. $p:G \to\bbr$, for every $ \a \in \Aut (G)$ and $a \in G$, $p (\a (a)) = a$ and $\sum\limits_{a\in G} p(a) = 1$)  there exists a word $w \in F_2$ such that for every $c\in G, \left | \frac{|\{ (a, b) \in G\times G | w(a,b) = c\} |}{|G|^2} - p(c)\right| < \vare$.

\medskip

\noindent {\bf 3)} \  The proof of the Theorem needs the classification only to ensure that every $a\in A' = A \backslash \{ e\}$ is part of a generating  2-set.  Without the classification the proof shows that if $A'$ is a set of elements such that every $a \in A'$ is part of a generating 2-set then $A = A' \cup \{ e\}$ is the image of some word map.  It follows that Corollary 2 does not rely on the classification. (We thank M. Kassabov for this last observation.)

\medskip

\noindent {\bf  4)} \  It is not so easy to find  a word $w$ as in Corollary 2 explicitly.  Our proof is effective in the sense that we can bound the length of $w$ in the proof.  For example, $|G|^{|G|^2}$ is a bound; (maybe by using some ideas and results on the uniformity of expanders one can do better) but these bounds are too large to be useful.  Our work suggests to study   quantitative versions of Theorem 1 and Corollary 2.   The work of Hadad [H] can be thought of as a step in this direction for $A = \{ e\}$.


\end{document}